\input amstex
\magnification=\magstep1
\comment\nologo
\endcomment
\document
\centerline{\bf{DUALITY BETWEEN PRIME FACTORS AND THE}}
\smallskip
\centerline{\bf{PRIME NUMBER THEOREM FOR ARITHMETIC PROGRESSIONS - II}}
\medskip
\centerline{\it{Krishnaswami Alladi and Jason Johnson}}
\medskip
\centerline{\it{Dedicated to George Andrews and Bruce Berndt for their 85th birthdays}}
\bigskip
ABSTRACT: {\it{In the first paper under this title (1977), the first author
utilized a duality identity between the largest and smallest prime factors
involving the Moebius function, to establish the following result as a consequence of the Prime Number Theorem for Arithmetic Progressions: If}} $k$ {\it{and}} $\ell$ {\it{are
positive integers, with}} $1\le\ell\le k$ {\it{and}} $(\ell, k)=1$, {\it{then}}
$$
\sum_{n\ge 2,\, p(n)\equiv\ell(mod\,k)}\frac{\mu(n)}{n}=\frac{-1}{\phi(k)},
$$
{\it{where}} $\mu(n)$ {\it{is the Moebius function,}} $p(n)$ {\it{is the smallest prime factor
of}} $n$, {\it{and}} $\phi(k)$ {\it{is the Euler function. Here we utilize the next level Duality identity between the second largest prime factor and the smallest prime factor, involving the Moebius function and}} $\omega(n)$, {\it{the number of distinct prime factors of}} $n$, {\it{to establish the following result as a consequence of the Prime Number Theorem for Arithmetic Progressions: For all}} $\ell$ {\it{and}} $k$ {\it{as above,}} 
$$
\sum_{n\ge 2, \, p(n)\equiv\ell(mod\,k)}\frac{\mu(n)\omega(n)}{n}=0.
$$
{\it{A quantitative version of this result is proved.}}

\medskip


\medskip


\bigskip

\centerline{\bf{\S1: Background}}
\medskip
Two famous results of Edmund Landau are that
$$
M(x):=\sum_{1\le n\le x}\mu(n)=o(x), \quad \text{as} \quad x\to\infty\tag1.1
$$
and
$$
\sum^{\infty}_{n=1}\frac{\mu(n)}{n}=0\tag1.2
$$
are (elementarily) equivalent to the Prime Number Theorem (PNT), where $\mu(n)$ is the
Moebius function. Similarly, there are results equivalent to the Prime Number Theorem for
Arithmetic Progressions (PNTAP) in which $\mu(n)$ is replaced by $\mu(n)\chi(n)$,
where $\chi(n)$ is a Dirichlet character modulo $k$, when the arithmetic progression under consideration has common difference $k$. 

In [2], the first author noticed the following interesting Duality identities
involving the Moebius function that connect the smallest and largest prime
factors of integers:

$$
\sum_{2\le d|n}\mu(d)f(p(d))=-f(P(n)),\tag1.3
$$
and
$$
\sum_{2\le d|n}\mu(d)f(P(d))=-f(p(n)),\tag1.4
$$
where for $d>1$, $p(d)$ and $P(d)$ denote the smallest and largest prime factors of $d$ respectively, and $f$ is ANY function on the
primes. Using (1.3) and properties of the Moebius function, it was shown in
[2] that if $f$ is a bounded function on the primes such that
$$
\lim_{x\to\infty}\frac{1}{x}\sum_{1\le n\le x}f(P(n))=c,\tag1.5
$$
then
$$
\sum^{\infty}_{n=2}\frac{\mu(n)f(p(n))}{n}=-c,\tag1.6
$$
and vice-versa. This is a surprising generalization of Landau's result (1.2).
To realize this is a generalization, rewrite (1.2) as
$$
\sum^{\infty}_{n=2}\frac{\mu(n)}{n}=-1.\tag1.7
$$
Then (1.7) follows from (1.5) and (1.6) by taking $f(p)=1$ for all primes $p$.

Next it was shown in [2] that the PNTAP implies that the sequence $P(n)$ of
largest prime factors is uniformly distributed in the reduced residue classes
modulo a positive integer $k$. So if $f$ is chosen to be the characteristic function of primes in an arithmetic progression $\ell(mod\,k)$, then for
such $f$, (1.5) holds with $c=1/\phi(k)$, and therefore
$$
\sum_{n\ge 2, \, p(n)\equiv\ell(mod\,k)}\frac{\mu(n)}{n}=\frac{-1}{\phi(k)},
\tag1.8
$$
for ALL positive integers $k$ and any $\ell$ satisfying $(\ell, k)=1$. 
This is even more surprising because it gives a way of slicing the convergent series in (1.7) into $\phi(k)$ subseries all converging to the same value! As far as we know, this is the first example of slicing convergent series into equal
valued subseries. In the last few years, (1.8) has received considerable attention and has been generalized in the setting of algebraic number theory (see [5],
[8], [10], [16], [21], and [22]).

In [2], the following more general duality identities were noted: For a positive integer $k$ let $P_k(n)$ and $p_k(n)$ denote the $k$-th largest and $k$-th smallest prime factors $n$ respectively (defined by strict inequalities), if $n$ has at least $k$ distinct prime factors. Also let $\omega(n)$ denote the number of distinct prime factors of $n$. Then
$$
\sum^*_{1<d|n}\mu(d)f(P_k(d))=(-1)^k\binom{\omega(n)-1}{k-1}f(p_1(n)),\tag1.9
$$
and
$$
\sum^*_{1<d|n}\mu(d)f(p_k(d))=(-1)^k\binom{\omega(n)-1}{k-1}f(P_1(n)),\tag1.10
$$
where the * over the summation means that if $n$ has fewer than $k$ distinct
prime factors, then the sum is zero. (NOTE: When $k=1$, we often write $p(n)$ and
$P(n)$ in place of $p_1(n)$ and $P_1(n)$, respectively.) From (1.9) and (1.10), it follows by
Moebius inversion that
$$
\sum_{1<d|n}\mu(d)\binom{\omega(d)-1}{k-1}f(P_1(d))=(-1)^kf(p_k(n)),\tag1.11
$$
and
$$
\sum_{1<d|n}\mu(d)\binom{\omega(d)-1}{k-1}f(p_1(d))=(-1)^kf(P_k(n)).\tag1.12
$$
In these identities, we adopt the convention that $f(P_k(n))=f(p_k(n))=0$ if
$\omega(n)<k$, that is, if $n$ has fewer than $k$ distinct prime factors.

In this paper we will discuss consequences of (1.12) in the case $k=2$, that is
the identity
$$
\sum_{1<d|n}\mu(d)(\omega(d)-1)f(p(d))=f(P_2(n))\tag1.13
$$
and its implications.

As a start, analogous to (1.2), we establish (see Theorem 4 of \S2)that
$$
\sum^{\infty}_{n=1}\frac{\mu(n)\omega(n)}{n}=0.\tag1.14
$$
Since $\omega(1)=0$, (1.14) gives
$$
\sum^{\infty}_{n=2}\frac{\mu(n)\omega(n)}{n}=0,\tag1.15
$$
unlike (1.7). Next we establish (see Theorem 7 of \S5),
that for each positive integer $k$, the sequence $P_2(n)$ of second
largest prime factors, is uniformly distributed in the reduced residue
classes $\ell (mod\,k)$. From this uniform distribution result, (1.15), and
by choosng $f$ to be the characteristic function for primes $\equiv\ell(mod\, k)$, we prove (see Theorem 10 of \S6) that
$$
\sum_{n\ge 2, p_1(n)\equiv\ell(mod\, k)}\frac{\mu(n)\omega(n)}{n}=0,\tag1.16
$$
and this is the main result of the paper.

To establish Theorems 7 and 10, several auxiliary results are proved. All the
theorems in this paper are proved in quantitative form.

{\bf{Notations and Conventions:}} In what follows, $c, c_1, c_2, \cdots$ are absolute positive constants whose values
will not concern us. The $<<$ and $O$ notations are equivalent and will be
used interchangeably as is convenient. We also adopt the convention that
$$
f(x)<< g(x) \quad \text{means} \quad |f(x)|<K|g(x)|,
$$
with $x$ ranging in some domain, and $K$ a positive constant. Implicit constants are absolute
unless otherwise indicated with a subscript. Although our results can be established
with uniformity by allowing the modulus $k$ to grow slowly as a function of $x$,
we only consider here an arbitrary but fixed modulus $k$. The alphabet $n$ whether used as the argument of a function, or in a summation, will always be a positive integer. Also, any time we have a sum over $p$, or have $p$ as an argument
of a function, it is to be understood that $p$ is prime.   

By $E(x,k,\ell)$ we mean the difference
$$
E(x,k,\ell)=\pi(x,k,\ell)-\frac{\ell i(x)}{\phi(k)},
$$
where
$$
\pi(x,k,\ell)
=\sum_{p\le x, p\equiv\ell(mod\,k}1, \quad \text{and} \quad
\ell i(x)=\int^x_2\frac{dt}{\log t}.
$$
When $k$ and $\ell$ are specific, we simply use $E(x)$ in place of
$E(x,k,\ell)$. Finally, by $R(x)$ we mean any decreasing function of $x$
that tends to zero as $x\to\infty$ and 
bounds from above the relative error in the PNTAP; that is
$$
|\pi(x,k,\ell)-\frac{\ell i(x)}{\phi(k)}|<\frac{\ell i(x)}{\phi(k)}R(x).
$$
In what follows, we will choose
$$
R(x)=e^{-c\sqrt log\,x},
$$
where $c$ can be any positive constant. Also $T=T(x)$ will be a function which
will be chosen optimally to get suitable
bounds in various estimates, but $T$ will not necessarily be the same in different contexts. We shall use the standard notation $[x]$ for the integral part
of a real number $x$, and $\{x\}$, where indicated, will denote the fractional part of $x$,
namely $x-[x]$. Finally, complex numbers will be denoted either by $z$, or
by $s=\sigma+it$ when dealing with Dirichlet series. 
Further notation will be introduced in the sequel as needed. 

\newpage

\centerline{\bf{\S2. The Moebius function and the number of prime factors}}

\medskip

Here and throughout, by $C=C(x)$ we mean the rectangular contour whose corners
in the complex plane are given by
$$
(1+\frac{1}{log\,x}, -T), \quad (1+\frac{1}{log\,x}, T), \quad
(-\frac{1}{log\,T}, T) \quad (-\frac{1}{log\,T}, -T),\tag2.1
$$
where $x\ge 3$ and $T=exp{\sqrt log\,x}$. Inside and on this contour, we know
that for some absolute constant $c_1>0$, 
$$
|\frac{1}{\zeta(s)}|<<(log\,x)^{c_1}, \tag2.2
$$
where $\zeta(s)$ is the Riemann zeta function. Since
$$
\sum^{\infty}_{n=1}\frac{\mu(n)}{n^s}=\frac{1}{\zeta(s)}, \quad \text{for} \quad
Re(s)>1,\tag2.3
$$
the Perron integral method applied to the contour $C$, together with (2.2)
gives
$$
M(x):=\sum_{n\le x}\mu(n)<<xe^{-c_2\sqrt log\,x}.\tag2.4
$$
Instead of summing over all positive integers $n$ as in (2.3), if we sum only
over those integers which are not multiples of a certain prime $p$, then we have
$$
\sum_{n\ge1,\, (n,p)=1}\frac{\mu(n)}{n^s}
=(1-\frac{1}{p^s})^{-1}\frac{1}{\zeta(s)}, \quad \text{for} \quad Re(s)>1.
\tag2.5
$$
Now on the contour $C$, we have the bound
$$
|(1-\frac{1}{p^s})^{-1}|\le (1-\frac{1}{p^{\sigma}})^{-1}\le
(1-\frac{1}{2^{\sigma}})^{-1}\le(1-\frac{1}{2^{\sqrt log\,3}})^{-1}\tag2.6
$$
valid uniformly for ALL primes $p$, where $\sigma=1-1/log\,T$. 
Thus the Perron integral method that yielded the bound for $M(x)$ in (2.4),
now gives in view of (2.6), the following bound
$$
M^{(p)}(x):=\sum_{n\le x, (n,p)=1}\mu(n)<< xe^{-c_2\sqrt log\,x}, \tag2.7
$$
valid uniformly for ALL primes $p$.

Next we consider the sum
$$
M_{\omega}(x):=\sum_{1\le n\le x}\mu(n)\omega(n).
$$
Even though our focus here is $n$ square-free, for the purpose of employing the  Moebius inversion formula, it is convenient to also consider the function
$\Omega(n)$, which is the (total) number of prime factors of $n$ counted with
multiplicity, because the function $\Omega$ is a totally additive. 

If $\chi_P$ denotes the characteristic function of the prime powers, then
$$
\Omega(n)=\sum_{d|n}\chi_P(d).
$$
Thus by Moebius inversion, we have
$$
\chi_P(n)=\sum_{d|n}\mu(d)\Omega(\frac{n}{d})
$$
$$
=\Omega(n)\sum_{d|n}\mu(d)-\sum_{d|n}\mu(d)\Omega(d)=-\sum_{d|n}\mu(d)\Omega(d),
\tag2.8
$$
because $\Omega(n)\sum_{d|n}\mu(d)$ is identically zero. Since $\mu(d)=0$
if $d$ is not square-free, we may rewrite (2.8) as
$$
\sum_{d|n}\mu(d)\omega(d)=-\chi_P(n),\tag2.9
$$
and Moebius inversion applied to (2.9) yields
$$
\mu(n)\omega(n)=-\sum_{d|n}\chi_P(d)\mu(\frac{n}{d}).\tag2.10
$$

Our first result is:
\medskip
{\bf{Theorem 1:}} {\it{For}} $x\ge 2$, {\it{we have}}
$$
M_{\omega}(x):=\sum_{n\le x}\mu(n)\omega(n) << \frac{x}{\log\, x}.
$$
\smallskip
{\bf{Proof:}} Use 
$$
M_{\omega}(x)=\sum_{n\le x}\mu(n)\omega(n)=\sum_{n\le x}\mu(n)\sum_{p|n}1
$$
$$
=\sum_{p\le x}\quad\sum_{n\le x, \, n\equiv 0(mod\, p)}\mu(n).\tag2.11
$$
In (2.11), in the inner sum on the right, put $n=mp$, with $(m,p)=1$, to
rewrite it as
$$
M_{\omega}(x)=-\sum_{p\le x}\quad\sum_{m\le x/p,\, (m,p)=1}\mu(m)=
-\sum_{p\le x}M^{(p)}(\frac{x}{p})=\Sigma_1+\Sigma_2,\tag2.12
$$
where
$$
\Sigma_1=\sum_{p\le T}M^{(p)}(\frac{x}{p})\quad \text{and}\quad
\Sigma_2=\sum_{T<p\le x}M^{(p)}(\frac{x}{p}),\tag2.13
$$ 
and $T$ will be chosen optimally below.

To estimate $\Sigma_1$ and $\Sigma_2$ we use the well-known estimate
$$
\sum_{p\le x}\frac{1}{p}=\log\log\,x+c_3+O(e^{-c_4\sqrt\log\,x}),\quad \text{for}
\quad x\ge 3,\tag2.14
$$
and its consequence
$$
\sum_{y<p\le x}\frac{1}{p}=\log\log\,x-\log\log\,y+O(e^{-c_4\sqrt\log\,y}),
\quad \text{for}\quad 3\le y\le x.\tag2.15
$$
So from (2.7) and (2.14) it follows that 
$$
\Sigma_1<<\sum_{p\le T}\frac{x}{p}e^{-c_2\sqrt\log(x/p)}
<<x\log\log\,x\,e^{-c_2\sqrt\log(x/T)}.\tag2.16
$$
To estimate $\Sigma_2$, we break it up as follows:
$$
\Sigma_2=\sum_{m<(x/T)-1}\quad\sum_{\frac{x}{m+1}<p\le \frac{x}{m}}, \tag2.17
$$
and denote the inner sum in (2.17) as $\Sigma^{(m)}$. 
We will set $T=x^{1-\varepsilon}$ and will choose $\varepsilon\to 0$ optimally
as $x\to\infty$. Thus we get the following bound for 
$\Sigma^{(m)}$ by using (2.7) and (2.15): 
$$
\Sigma^{(m)}<<xe^{-c_2\sqrt\log\,m}\{(\log\log\frac{x}{m}-\log\log\frac{x}{m+1})+O(e^{-c_4\sqrt\log(x/m)})\}
$$
$$
<<xe^{-c_2\sqrt\log\,m}\frac{1}{m\log\,x}+xe^{-c_2\sqrt\log\,m}e^{-c_4\sqrt\log(x/m)}.\tag2.18
$$
Note that $\frac{x}{T}=x^{\varepsilon}$, and that
$$
\sum^{\infty}_{m=1}\frac{e^{-c_2\sqrt\log\,m}}{m}<\infty.
$$
Thus by summing the expression on the right in (2.18) over $m\le x^{\varepsilon}$, we get
$$
\Sigma_2<<\frac{x}{\log\,x}+x\sum_{m<x^{\varepsilon}}e^{-c_2\sqrt\log\,m}e^{-c_4\sqrt\log(x/m)}
$$
$$
<<\frac{x}{\log\,x}+x^{1+\varepsilon}e^{-\sqrt(1-\varepsilon)\log\,x}.\tag2.19
$$
Finally, we choose 
$$
\varepsilon=c_5\frac{(\log\log\,x)^2}{\log\,x},\quad \text{with} \quad
c_5=4c^{-2}_2.\tag2.20
$$
So (2.16), (2.19) and (2.20) yield
$$
\Sigma_1<<\frac{x\log\log\,x}{\log^2x} \quad \text{and} \quad \Sigma_2<<\frac{x}{\log\,x}.\tag2.21
$$
Theorem 1 follows from (2.12) and (2.21). 

\medskip

{\bf{Remark:}} If we choose $c_5=N^2c^{-2}_2$, with $N$ arbitrarily
large, then we would get
$$
\Sigma_1<<_N\frac{x\log\log\,x}{\log^Nx},
$$
but this is of no use since $\Sigma_2$ is
bounded only by $x/\log\,x$, and not any better by the above method.
For a sharper estimate for $M_{\omega}(x)$ due to Tenenbaum by analytic methods
see (3.15) below.  

Next, using Theorem 1 and following an idea of Axer, we prove:

\medskip

{\bf{Theorem 2:}} {\it{With}} $\{w\}$ {\it{denoting the fractional part of}} $w$, {\it{we have}}  
$$
\sum_{n\le x}\mu(n)\omega(n)\{\frac{x}{n}\}<<
\frac{x\sqrt\log\log\,x}{\sqrt\log\,x}.
$$

{\bf{Proof:}} Begin by splitting
$$
\sum_{n\le x}\mu(n)\omega(n)\{\frac{x}{n}\}=\sum_{n\le T}+\sum_{T<n\le x}
=\Sigma_3+\Sigma_4,\tag2.22
$$
with $T=T(x)$ to be determined below.

For $\Sigma_3$, we use the trivial bound
$$
|\Sigma_3|\le \sum_{n\le T}\omega(n)<< T\log\log\,T.\tag2.23
$$
To estimate $\Sigma_4$, we use partial summation:
$$\Sigma_4=\sum_{T<n\le x}(M_{\omega}(n)-M_{\omega}(n-1))\{\frac{x}{n}\}
$$
$$
=\sum_{T<n\le {x-1}}M_{\omega}(n)(\{\frac{x}{n}\}-\{\frac{x}{n+1}\})
+O(|M_{\omega}(x)|)+(|M_{\omega}(\frac{x}{T})|).\tag2.24
$$
From (2.24) and Theorem 1 we deduce that
$$
|\Sigma_4|<< \frac{x}{log\,x}+\sum_{T<n\le {x-1}}|M_{\omega}(n)||\{\frac{x}{n}\}-\{\frac{x}{n+1}\}|
<<\frac{x}{\log\,x}\sum_{T<n\le {x-1}}|\{\frac{x}{n}\}-\{\frac{x}{n+1}\}|.\tag2.25
$$
At this stage, we note that
$$
\sum_{T<n\le x}|\{\frac{x}{n}\}-\{\frac{x}{n+1}\}|\le V_{\{\}}[1,\frac{x}{T}]
<<\frac{x}{T},\tag2.26
$$
where $V_{\{\}}[a,b]$ is the total variation of $\{x\}$ on $[a,b]$. So by (2.25) and (2.26) we have
$$
\Sigma_4<<\frac{x}{\log\,x}\frac{x}{T}.\tag2.27
$$
Thus (2.22), (2.23) and (2.27) yield
$$
\sum_{n\le x}\mu(n)\omega(n)\{\frac{x}{n}\}<<
T\log\log\,T+\frac{x}{\log\,x}\frac{x}{T}.\tag2.28
$$
We need to choose $T$ optimally to minimize the right hand side of (2.28). 
The choice
$$
T=\frac{x}{\sqrt{\log\,x\,\log\log\,x}}
$$
in (2.28) yields
$$
\sum_{n\le x}\mu(n)\omega(n)\{\frac{x}{n}\}<<
\frac{x\sqrt\log\log\,x}{\sqrt\log\,x}
$$
which proves Theorem 2. 

Since Theorem 2 deals with the fractional part function as the weight, we establish next the corresponding result with the weight as the integral part function:

\medskip

{\bf{Theorem 3:}}
$$
\sum_{n\le x}\mu(n)\omega(n)[\frac{x}{n}]=\frac{x}{\log\,x}+O(\frac{x}{\log^2x}).
$$

{\bf{Proof:}} Note that (2.9) yields
$$
\sum_{n\le x}\mu(n)\omega(n)[\frac{x}{n}]=\sum_{n\le x}\sum_{d|n}\mu(d)\omega(d)
$$
$$
=\sum_{n\le x}\chi_P(n)=\frac{x}{\log\,x}+O(\frac{x}{\log^2x}),
$$
which proves Theorem 3.

Theorems 2 and 3 lead to the the main result of this section:

\medskip

{\bf{Theorem 4:}}
$$
m_{\omega}(x):= \sum_{n\le x}\frac{\mu(n)\omega(n)}{n}=O(\frac{\sqrt\log\log\,x}{\sqrt\log\,x}).
\tag2.29
$$
In particular
$$
\sum^{\infty}_{n=1}\frac{\mu(n)\omega(n)}{n}
=\sum^{\infty}_{n=2}\frac{\mu(n)\omega(n)}{n}=0.
$$

{\bf{Proof:}} From Theorems 2 and 3, we get
$$
\sum_{n\le x}\mu(n)\omega(n)\frac{x}{n}=\sum_{n\le x}\mu(n)\omega(n)[\frac{x}{n}]
+\sum_{n\le x}\mu(n)\omega(n)\{\frac{x}{n}\}
$$
$$
<<\frac{x\sqrt\log\log\,x}{\sqrt\log\,x}+\frac{x}{\log\,x}.\tag2.30
$$
By cancelling $x$ on both extremes of (2.30), we get (2.29) of Theorem 4. By
letting $x\to\infty$ in (2.29), we get (1.14) and (1.15), which is the second
assertion of Theorem 4.

{\bf{Remark:}} When I communicated the results of this section to Tenenbaum, he responded [19] by saying that he can establish stronger quantitative versions of the Theorems 1 and 4 by the Selberg-Delange analytic method. We describe Tenenbaum's
approach and state his stronger quantitative results in the next section. But we mention here that the actual results of $\S2$ and $\S3$ are not used to establish our main result in $\S6$, but the elementary method used in this section is what is employed in $\S6$. The discussion on $M_{\omega}(x)$ and $m_{\omega}(x)$ in
this and the next section provides a context to understand our main result. 

\bigskip

\centerline{\bf{\S3. Analytic approach to sums of $\mu(n)\omega(n)$}}

\medskip

For a complex number $z$, represent the function $\zeta(s)^z$ as a Dirichlet series
$$
\zeta(s)^z=\sum^{\infty}_{n=1}\frac{d_z(n)}{n^s}, \quad \text{for} \quad Re(s)>1.
\tag3.1
$$
The function $d_z(n)$ is called the generalized divisor function because
$d_2(n)=d(n)$ is the standard divisor function. In a fundamental paper,
Selberg [14] showed that
$$
\sum_{n\le x}d_z(n)=\frac{x(\log\,x)^{z-1}}{\Gamma(z)}+O_R(x(\log\,x)^{z-2})\tag3.2
$$
is valid uniformly for $|z|\le R$.

In order to establish (3.2), Selberg used the Perron integral method, but
since $\zeta(s)^z$ has a branch point singularity at $s=1$ when $z$ is not an
integer, he modified the contour $C(x)$ in \S2 by replacing the short
line segment from $b-i(2\log\,T)^{-1}$ to $b+i(\log\,T)^{-1}$ on $C(x)$ with
the following lacet $L$ around $s=1$:

$L$ starts at $b-i(2\log\,x)^{-1}$, runs parallel to the $x$-axis until
$1-i(2\log\,x)^{-1}$, then encircles $s=1$ in a semi-circle of radius
$(2\log\,x)^{-1}$, and ends with a line segment from $1+i(\log\,x)^{-1}$
to $b+i(\log\,x)^{-1}$. 

The contribution around $L$ leads to the Hankel contour for the Gamma function,
and that explains the presence of $\Gamma(z)$ in (3.2).

Selberg's method applies more generally to sums of coefficients of Dirichlet
series representing functions of the type
$$
\zeta(s)^z. H(s),
$$
where $H(s)$ would be analytic in the half plane $Re(s)>\frac{1}{2}$, and
indeed Selberg considered a few important such $H(s)$ in [14]. A case of interest
to us here is the sum
$$
S_{-z}(x):=\sum_{n\le x}\mu(n)z^{\omega(n)}\tag3.3
$$
which can be viewed as the sum of the coefficients of the Dirichlet series
$$
\sum^{\infty}_{n=1}\frac{\mu(n)z^{\omega(n)}}{n^s}=\prod_p(1-\frac{z}{p^s}) \quad
\text{for} \quad Re(s)>1.\tag3.4
$$
We could rewrite (3.4) as
$$
\sum^{\infty}_{n=1}\frac{\mu(n)z^{\omega(n)}}{n^s}=\zeta(s)^{-z}G(s,z),\tag3.5
$$
where
$$
G(s,z)=\prod_p(1-\frac{z}{p^s})(1-\frac{1}{p^s})^{-z},\tag3.6
$$
is analytic in $Re(s)>\frac{1}{2}$, and uniformly bounded in the half plane
$Re(s)\ge \frac{1}{2}+\delta$, for each $\delta>0$. So by the Selberg method,
one gets
$$
S_{-z}(x)=\frac{x(\log\,x)^{-z-1}G(1,z)}{\Gamma(-z)}
+O_R(x(\log\,x)^{-z-2}),\tag3.7
$$
is uniformly valid for $|z|\le R$, and indeed this is implicit in [13]. 

Selberg's method was extended by Delange to deal with sums of coefficients
of Dirichlet series convergent in $Re(s)>1$, and can be represented as
$$
\zeta(s)^z(\log\zeta(s))^kH(s),\tag3.8
$$
where $k$ is a non-negative integer, and $H(s)$ is analytic in
$Re(s)>\frac{1}{2}$. This is a natural extension of Selberg's method (and one
to be expected) because
$\zeta^z(s)=e^{z\log\zeta(s)}$, but it is a useful extension. Tenenbaum [18] has
a thorough account of the Selberg-Delange method in its most general form, and in doing so, has improved on the quantitative aspects as well. 

We can view $M_{\omega}(x)$ as
$$
M_{\omega}(x)=\frac{d}{dz}S_{-z}(x)|_{z=1}.\tag3.9
$$
So one may heuristically get from (3.7) that
$$
M_{\omega}(x)\sim \frac{d}{dz}(\frac{x(\log\,x)^{-z-1}G(1,z)}{\Gamma(-z)})|_{z=1}.
\tag3.10
$$
Note that, 
$$
\frac{G(1,z)}{\Gamma(-z)}
$$
has a simple zero at $z=1$, and so (3.10) would yield 
$$
M_{\omega}(x)\sim \frac{c_7x}{(\log\,x)^2},\tag3.11
$$
with some non-zero constant.
But Tenenbaum in his letter [19] established (3.11) by the Selberg-Delange
method, as well as a precise series expansion as detailed below. 

Observe that
$$
\sum^{\infty}_{n=1}\frac{\mu(n)\omega(n)}{n^s}=
\frac{d}{dz}(\sum^{\infty}_{n=1}\frac{\mu(n)z^{\omega(n)}}{n^s})|_{z=1},
\quad \text{for} \quad Re(s)>1,\tag3.12
$$
because the term-by-term differentiation of the Dirichlet series on the right
is valid in $Re(s)>1$. So by (3.5) and (3.12) we have
$$
\sum^{\infty}_{n=1}\frac{\mu(n)\omega(n)}{n^s}=
\frac{d}{dz}(\zeta^{-z}(s)G(s,z))|_{z=1}
$$
$$
=-\zeta^{-z}(s)\log\zeta(s)G(s,z)|_{z=1}+\zeta^{-z}(s)G'(s,z)|_{z=1}
$$
$$
= -\zeta^{-1}(s)\log\zeta(s)G(s,1)+\zeta^{-1}(s)G'(s,1), \quad \text{for} \quad
Re(s)>1.\tag3.13
$$
In view of the representation in (3.13) for the Dirichlet series, the Selberg-Delange method can be applied using the function on the right in (3.13), to deduce that
$$
M_{\omega}(x)= \frac{c_7x}{\log^2x}+O(\frac{x}{\log^3x}).\tag3.14
$$
Tenenbaum [19] notes that Theorem 11.5.2 in his book [18] readily yields the
following more precise estimate:
$$
M_{\omega}(x)=x\sum_{0\le k\le N}\frac{\lambda_k}{\log^{k+2}x}+O(xR_{N+2}(x)),
\tag3.15
$$
where the $\lambda_k$ are constants, with $\lambda_0=1$ (consequently $c_7=1$, and 
$$
R_N(x)=e^{-\sqrt\log\,x}+O((\frac{c_8N+1}{\log\,x})^{N+1}).\tag3.16
$$
But, as noted in the previous section, such a precise result as (3.16) for $M_{\omega}(x)$ is not needed for us here. 

It follows from (3.15) by partial summation that
$$
\sum_{n\le x}\frac{\mu(n)\omega(n)}{n}= c_9+O(\frac{1}{\log\,x}),
$$
and so the series
$$
\sum^{\infty}_{n=1} \frac{\mu(n)\omega(n)}{n}
$$
is convergent to $c_9$. If we call the expression on the right in (3.13) as $F(s,1)$, then Tenenbaum [19] notes that
$$
c_9=\lim_{\sigma\to 1^+}F(\sigma,1)=0
$$
and so
$$
\sum_{n\le x} \frac{\mu(n)\omega(n)}{n}<<\frac{1}{\log\,x}\tag3.17
$$
which is stronger than our Theorem 4, which we derived elementarily from the
strong form of the Prime Number Theorem.

Indeed, in the spirit of (3.15), Tenenbaum [19] proved the stronger result
$$
\sum_{n\le x}\frac{\mu(n)\omega(n)}{n}=\sum_{1\le k\le N}\frac{\nu_k}{\log^kx}
+O(R_{N+1}(x)).\tag3.18
$$
which follows by directly applying the Selberg-Delange method to evaluate
the sum by considering the associated Dirichlet series
$$
\sum^{\infty}_{n=1}\frac{\mu(n)z^{\omega(n)}}{n^{s+1}}.
$$

\newpage

{\bf{Some key differences:}}

Denote by
$$
m(x):=\sum_{n\le x}\frac{\mu(n)}{n}.\tag3.19
$$
It is known by a theorem of Landau on Dirichlet
series whose coefficients are eventually of the same sign (with a similar
result for Dirichlet type integrals), that $M(x)$ changes sign infinitely
often, and that $m(x)$ changes sign infinitely often as it converges to 0
as $x\to\infty$. In contrast, in view of the fact that $M_{\omega}(x)$ can
be estimated asymptotically with a leading term as in (3.14), $M_{\omega}(x)$
will eventually be of the same sign, and therefore will NOT change
sign infinitely often. Similarly, 
$$
m_{\omega}(x):=\sum_{n\le x}\frac{\mu(n)\omega(n)}{n}\tag3.20
$$
will be eventually of the same sign as it tends to 0 when $x\to\infty$, and
so will NOT change sign infinitely often.

The other key difference is in their sizes. Whereas the strong form of the
Prime Number Theorem implies that
$$
M(x)=O(xe^{-c\sqrt\log\,x}) \quad \text{and} \quad
m(x)=O(e^{-c\sqrt\log\,x}),
$$
we have
$$
M_{\omega}(x)\sim c_7\frac{x}{\log^2x} \quad \text{and} \quad
m_{\omega}(x)\sim c_{10}\frac{1}{\log\,x}.
$$
It is however to be noted that by writing
$$
\sum^{\infty}_{n=1}\frac{\mu(n)\omega(n)}{n}=
\int^{\infty}_{1^-}\frac{dM_{\omega}(x)}{x},\tag3.21 
$$
and using integration-by-parts, we get
$$
\sum^{\infty}_{n=1}\frac{\mu(n)\omega(n)}{n}=
\int^{\infty}_1\frac{M_{\omega}(x)}{x^2}dx=0.\tag3.22
$$
Thus $M_{\omega}(x)$ must change sign, even though it does not change sign infinitely often!

\bigskip

\centerline{\bf{\S4: The sizes of the largest and second largest prime factors}}
\medskip
The fundamental counting function associated with the largest prime factor
$P(n)$ is
$$
\Psi(x,y)=\sum_{n\le x. P(n)\le y}1.\tag4.1
$$
Here and in what follows, we shall denote by $\alpha$, the quantity
$\frac{\log\,x}{\log\,y}$. In an important paper [7], de Bruijn showed that
with some constant $c>0$, 
$$
\Psi(x,y)<< xe^{-c\alpha}, \quad \text{uniformly for} \quad 2\le y\le x.
\tag4.2a
$$
Tenenbaum [19; Theorem III.5.1] has shown that (4.2a) holds with $c=1/2$. de Bruijn also proved that
$$
\Psi(x,y)<<x\log^2ye^{-\alpha\log\,\alpha-\alpha\log\log\,\alpha+O(\alpha)}, \quad  \text{for} \quad y>\log^2x,
\tag4.2b
$$
and indeed the uniform asymptotic estimate
$$
\Psi(x,y)\sim x\rho(\alpha),\quad \text{for} \quad e^{(\log\,x)^{3/5}}\le y\le x,
\tag4.2c
$$
where $\rho$ satisfies the integro-difference equation
$$
\rho(\alpha)=1-\int^{\alpha}_1\frac{\rho(u-1) du}{u}\tag4.3a
$$
and (de Bruijn [6])
$$
\rho(\alpha)=e^{-\alpha\log\,\alpha - \alpha\log\log\,\alpha + O(\alpha)}.\tag4.3b
$$
Thus $\Psi(x,y)$ is quite small in comparison with $x$ when $\alpha$ is large.

{\bf{Remark:}} Note that (4.2b) is of no use when $\alpha >1$ is fixed, because trivially $\Psi(x,y)\le x$ (!); so (4.2b) is used only when $\alpha\to\infty$ with $x$.
Much better extimates for $\Psi(x,y)$ are known including those that significantly extend the range of the asymptotic formula (4.2c); for such
superior results, see Hildebrand and Tenenbaum [9]. For our purpose
here, these superior results on $\Psi(x,y)$ are not needed; the above bounds suffice.

Next consider $P_2(n)$, the second largest prime factor of $n$. Note that
whereas $P(n)$ is uniquely defined, there are two ways to define the second
largest prime factor. We could define $P_2(n)=P(n/P(n))$ or $P_2(n)$ as the
largest prime factor of $n$ strictly less than $P(n)$. In the former definition, we
set $P_2(n)=1$ if $\Omega(n)<2$, and in the latter definition we set $P_2(n)=1$,
if $\omega(n)<2$. From the point of view of asymptotic estimates, there is
little difference between the two ways of defining $P_2(n)$. This is made
precise by:

\medskip

{\bf{Theorem 5:}} {\it{Let}} $N(x)$ {\it{denote the number of positive integers}} $n\le x$ {\it{for which}} $P(n)$ {\it{repeats. Then}}
$$
N(x)<<\frac{x}{e^{(\frac{1}{2}+o(1))\sqrt(\log\,x\log\log\,x)}}.
$$

{\bf{Proof:}} By (4.2b),
$$
\Psi(x,e^{\sqrt(\log\,x\log\log\,x})<< \frac{x}{e^{(\frac{1}{2}+o(1))\sqrt((\log\,x\log\log\,x)}}.\tag4.4
$$
So it suffices to consider those integers $n$ for which
$P(n)>e^{\sqrt(\log\,x\log\log\,x)}$.
Among these integers $\le x$, the number of those with largest prime $P(n)=p$
repeating, is trivially $O(x/{p^2})$. Thus
$$
N(x)<<\Psi(x,e^{\sqrt(\log\,x\log\log\,x)})+\sum_{p>exp{\sqrt(\log\,x}\log\log\,x)}\frac{x}{p^2}
$$
$$
<< \frac{x}{e^{(\frac{1}{2}+o(1))\sqrt(\log\,x\log\log\,x)}}+\sum_{n>exp{\sqrt(\log\,x\log\log\,x)}}\frac{x}{n^2}
<<\frac{x}{e^{(\frac{1}{2}+o(1))\sqrt(\log\,x\log\log\,x)}},
$$
which proves Theorem 5. 

Here we shall use the definition for $P_2(n)$ as the largest prime factor strictly less than the largest prime factor. From
Theorem 5 we see that we can focus on those integers for which $P(n)$ occurs
square-free.

Consider the counting function
$$
\Psi_2(x,y)=\sum_{n\le x,\, P_2(n)\le y}1.\tag4.5
$$
In contrast to $\Psi(x,y)$ which is very small in comparison with $x$ when $\alpha$ is large, the function $\Psi_2(x,y)$ is not that small. To realize this, observe that all integers of the form $2p\le x$ where $p$ is prime, will have $P_2(n)=2$, and so
$$
\Psi_2(x,y)\ge \Psi_2(x,2)>>\frac{x}{\log\,x}, \quad \text{for all} \quad y\ge 2.\tag4.6
$$
What we need here is a quantitative version of the fact that 
for ``almost all'' integers, $P_2(n)$ is large. This (and much more) is provided by a result of Tenenbaum [], on the size of $P_k(n)$, when $k\ge 2$. We state
Tenenbaum's result for $k=2$  (his eqns (1.5) and (1.6) in []) in the form of
\medskip
{\bf{Theorem 6:}} (Tenenbaum) {\it{There exists a function}} $\rho_2(\alpha)$
{\it{such that}} 
$$
\Psi_2(x,y)=x\rho_2(\alpha)(1+O(\frac{1}{\log\,y})), \quad \text{uniformly for}
\quad 2\le y \le x.\tag4.7
$$
{\it{The function}} $\rho_2(\alpha)$ {\it{satisfies}}
$$
\frac{1}{\alpha}<<\rho_2(\alpha)<<\frac{1}{\alpha}.\tag4.8
$$

Tenenbaum's proof of his stronger quantitative result on the joint
distribution of the $P_k(n)$ for $k\ge 2$, is quite intricate, and makes use
of sharp estimates for $\Psi(x,y)$. He has a precise formula for $\rho_2(\alpha)$ which we do not need here. For our purpose, all we need is   

\medskip

{\bf{Theorem 6*:}} {\it{Uniformly for}} $2\le T\le x$, {\it{we have}}
$$
\Psi_2(x,T)<<\frac{x\log\,T}{\log\,x}.
$$
which follows from Theorem 6. But we point out that the bound in Theorem 6*
for $T\le exp\{(\log\,x)^{1-\delta}\}$, for any $\delta >0$, can be proved using
the bounds in (4.2a) and (4.2b); the implicit constant will depend on
$\delta$. Theorem 6* will be used in what follows. 

\bigskip

\centerline{\bf{\S5. The uniform distribution of $P_2(n)$ modulo $k$}}
\medskip
In this section, we shall prove:

\medskip

{\bf{Theorem 7}} {\it{For each integer}} $k\ge 2$, {\it{the sequence}} $P_2(n)$ {\it{of the second largest prime factors, is uniformly distributed in the reduced residue classes modulo k. More precisely, for each fixed}} $k\ge 2$,
{\it{and any}} $1\le\ell<k$ {\it{with}} $(\ell, k)=1$, {\it{we have}}
$$
N_2(x,k,\ell):=\sum_{n\le x,\, P_2(n)\equiv \,\ell(mod\ k)}1=\frac{x}{\phi(k)}+
O(\frac{x(\log\log\,x)^2}{\log\,x}).\tag5.1
$$
\medskip
{\bf{Remark:}} Note that the number of positive integers up to $x$ with $\omega(n)=1$ or $\Omega(n)=1$ is
$$
\pi(x)+O(\sqrt x)=\frac{x}{\log\,x}+O(\frac{x}{\log^2x})=o(x)\tag5.2
$$
Thus it does not matter in Theorem 7 whether the sum in (5.1) is taken over
all integers $n\le x$ for which $P_2(n)\equiv\ell (mod\,k)$, or restricted
to integers for which $\omega(n)\ge 2$. 

\medskip

{\bf{Proof:}} For a prime $p$, denote by $S_2(x,p)$ the set of integers
$n\le x$ for which $P_2(n)=p$. Then
$$
\sum_{p\le \sqrt x}|S_2(x,p)|=x- \frac{x}{\log\,x}+O(\frac{x}{\log^2x}).\tag5.3
$$
Now if $N\in S_2(x,p)$, then we may write
$$
N=mpq, \quad \text{where}  \quad q\ge p \quad \text{and} \quad P(m)\le p,
$$
with $q$ being prime. In particular $m\le (x/p^2)$. Thus
$$
|S_2(x,p)|=\sum_{m\le (x/p^2), \,P(m)\le p}\quad\sum_{q>p, \, mpq\le x}1=
\sum_{p<q\le x/p}\quad\sum_{m\le x/(pq),\, P(m)\le p}1.\tag5.4
$$
Thus by summing the expression in (5.4) over $p\le\sqrt x$, we get
$$
\sum_{p\le\sqrt x}|S_2(x,p)|=
\sum_{p\le\sqrt x}\,\,\sum_{p\le q\le x/p}\Psi(\frac{x}{pq}, p)
=\sum_{q\le (x/2)}\,\,\sum^*\Psi(\frac{x}{pq}, p),\tag5.5
$$
where the * over the inner sum on the right means that the conditions on
$p$ and $q$ are
$$
p\le\sqrt x,\quad p\le q, \quad \text{and} \quad p\le \frac{x}{q}.\tag5.6
$$
To deal efficiently with the double sum on the right in (5.5), we consider
two cases, namely $q\le \sqrt x$ and $q>\sqrt x$, together with the inequalities in (5.6). This gives
$$
\sum_{p\le\sqrt x}|S_2(x,p)|=
\sum_{q\le\sqrt x}\,\,\sum_{p\le q}\Psi(\frac{x}{pq},p)
+\sum_{\sqrt x<q\le (x/2)}\,\,\sum_{p\le (x/q)}\Psi(\frac{x}{pq},p).\tag5.7
$$
At this point we note that we have already established an effective version of the statement that $P_2(n)$ is ``almost always'' large, namely, Theorem 6*.
In view of Theorem 6*, we may consider in (5,7) only the integers $n\le x$ for
which $P_2(n)>y$, with $y$ to be chosen later to satisfy $\log \,y=o(\log\,x)$.
Thus we modify (5.7) to

$$
\sum:=\sum_{y<p\le\sqrt x}|S_2(x,p)|=
\sum_{y<q\le\sqrt x}\,\,\sum_{y<p<q}\Psi(\frac{x}{pq},p)
+\sum_{\sqrt x<q\le (x/y)}\,\,\sum_{y<p\le (x/q)}\Psi(\frac{x}{pq},p)
$$
$$
=x+O(\frac{x\log\,y}{\log\,x}).\tag5.8
$$

Next we shall compare $\sum$ with
$$
I:=
\sum_{y<q\le\sqrt x}\,\int^q_y\Psi(\frac{x}{tq},t)\frac{dt}{\log\,t}
+\sum_{\sqrt x<q\le (x/y)}\,\,\int^{x/q}_y\Psi(\frac{x}{tq},t)\frac{dt}{\log\,t}
\tag5.9
$$
and estimate the difference (error) $E=\sum \, - \,I$ by using the strong form of
the Prime Number Theorem. We first consider the absolute value of the difference
$$
E_1:=|\sum_{y<q\le\sqrt x}\{\sum_{y<p\le q}\Psi(\frac{x}{pq}, p)
-\int^q_y\Psi(\frac{x}{tq},t)\frac{dt}{\log\,t}\}|
$$
$$
=|\sum_{y<q\le\sqrt x}\{\sum_{y<p<q}\,\,\sum_{n\le (x/pq), P(n)\le p}1
-\int^q_y(\sum_{n\le (x/tq), P(n)\le t}1)\frac{dt}{\log\,t}\}|
$$
$$
\le \sum_{y<q\le\sqrt x}\,\,\sum_{n\le (x/yq)}|\sum_{max(P(n),y)\le p\le min(x/nq,q)}1
-\int^{min(x/nq,q)}_{max(P(n),y)}\frac{dt}{\log\,t}|.\tag5.10
$$
It is to be noted that for the final term in (5.10),
we have
dropped the condition on $P(n)$ in the sum over $n$
on the right for simplicity since this will not lessen the effectiveness of
the upper bound on $E_1$ that we will get. We now use the strong form of the Prime Number Theorem on the expression on the right in (5.10) to deduce that
$$
E_1<<\sum_{y<q\le\sqrt x}\,\,\sum_{n\le (x/yq)}\frac{x}{nq\,exp{\{\sqrt\log\,(x/nq)\}}}
\le \, \frac{x}{exp{(\sqrt\log\,y)}}\sum_{y<q\le\sqrt x}\frac{1}{q}
\,\sum_{n\le (x/yq)}\frac{1}{n}
$$
$$
<<\frac{x\log\,x}{exp{(\sqrt\log\,y)}}\sum_{y<q\le\sqrt x}\frac{1}{q}
<<\frac{x\log\,x\log\log\,x}{exp{(\sqrt\log\,y)}}.\tag5.11
$$
In obtaining this upper bound for $E_1$, we have used the fact that
$z/exp{\sqrt\log\,z}$ is an increasing function of $z$, and so in deriving (5.11) we only used the error term in the strong form of the Prime Number Theorem
with $z=x/nq$.

Similarly, we bound the difference

$$
E_2:=|\sum_{\sqrt x<q\le (x/y)}\{\sum_{y<p<(x/q)}\Psi(\frac{x}{pq}, p)
-\int^{x/q}_y\Psi(\frac{x}{tq},t)\frac{dt}{\log\,t}\}|
$$
$$
=|\sum_{\sqrt x<q\le (x/y)}\{\sum_{y<p<(x/q)}\,\,\sum_{n\le (x/pq), P(n)\le p}1
-\int^{x/q}_y(\sum_{n\le (x/tq), P(n)\le t}1)\frac{dt}{\log\,t}\}|
$$
$$
\le \sum_{\sqrt x<q\le (x/y)}\,\,\sum_{n\le (x/yq)}|\sum_{max(P(n),y)\le p\le (x/nq)}1
-\int^{x/nq}_{max(P(n),y)}\frac{dt}{\log\,t}|
$$
$$
<< \sum_{\sqrt x<q\le (x/y)}\,\,
\sum_{n\le (x/yq)}\frac{x}{nq\,exp{\{\sqrt\log\,(x/nq)\}}}
\le \, \frac{x}{exp{(\sqrt\log\,y)}}\sum_{\sqrt x<q\le(x/y)}\frac{1}{q}
\,\sum_{n\le (x/yq)}\frac{1}{n}
$$
$$
<<\frac{x\log\,x}{exp{(\sqrt\log\,y)}}\sum_{\sqrt x<q\le (x/y)}\frac{1}{q}
<<\frac{x\log\,x\log\log\,x}{exp{(\sqrt\log\,y)}}.\tag5.12
$$
So from (5.11) and (5.12), we see that
$$
|E|=|\sum \,-\,I|\le E_1+E_2<<
\frac{x\log\,x\log\log\,x}{exp{(\sqrt\log\,y)}}.\tag5.13
$$
From (5.13) and (5.8) we deduce that
$$
I=x+O(\frac{x\log\,y}{\log\,x})+
O(\frac{x\log\,x\log\log\,x}{exp{(\sqrt\log\,y)}}).\tag5.14
$$
At this stage, we make the choice
$$
y=exp{\{(2\log\log\,x)^2\}},\tag5.15
$$
to conclude that
$$
I=x+O(\frac{x(\log\log\,x)^2}{\log\,x}) \quad \text{and} \quad
\sum =x+O(\frac{x(\log\log\,x)^2}{\log\,x}).\tag5.16
$$
This will be crucial in establishing Theorem 7.  

Now for an arbitrary but fixed integer $k\ge 2$, and for any $1\le \ell<k$ with
$(\ell, k)=1$, we consider the set $S^{k,\ell}_2(x)$ of integers $n\le x$ such
that, $\omega(n)\ge 2$, and 
$P_2(n)\equiv \ell(mod\,k)$. By classifying the members of this set in terms
of their second largest prime factor, we see that
$$
N_2(x,k,\ell)= |S^{k,\ell}_2(x)|=\sum_{p\le\sqrt x,\, p\equiv\ell(mod\,k)}|S_2(x,p)|.\tag5.17.
$$
In view of Theorem 6*, we have for any $y\le exp{(\log\,x)^{1-\delta}}$
$$
\sum_{p\le y, \, p\equiv\ell(mod\,k)}|S_2(x,p)|\le \sum_{p\le y}|S_2(x,p)|<<
\frac{x\log\,y}{\log\,x}.\tag5.18
$$
Thus
$$
N_2(x,k,\ell)=|S^{k,\ell}_2(x)|=\sum_{y<p\le\sqrt x,\,  p\equiv\ell(mod\,k)}|S_2(x,p)|
+O(\frac{x\log\,y}{\log\,x})\tag5.19
$$

Let us denote the sum on the right of (5.19) as $\sum^{k,\ell}$. Then by reasoning as above, we get
$$
\sum^{k,\ell}=
\sum_{y<q\le\sqrt x}\,\,\sum_{y<p<q,\, p\equiv\ell(mod\,k)}\Psi(\frac{x}{pq},p)
+\sum_{\sqrt x<q\le (x/y)}\,\,\sum_{y<p\le (x/q),\, p\equiv\ell(mod\,k)}\Psi(\frac{x}{pq},p).\tag5.20
$$
We now want to compare the expression in (5.20) with
$$
I^{k,\ell}:=
\sum_{y<q\le\sqrt x}\,\,\frac{1}{\phi(k)}\int^q_y\Psi(\frac{x}{tq},t)\frac{dt}{\log\,t}
+\sum_{\sqrt x<q\le (x/y)}\,\,\frac{1}{\phi(k)}\int^{x/q}_y\Psi(\frac{x}{tq},t)\frac{dt}{\log\,t}
\tag5.21
$$
Clearly
$$
I^{k, \ell}=\frac{I}{\phi(k)},\tag5.22
$$
with $I$ as in (5.9). Just as we obtained the bound in (5.13) for the difference
$\sum -I$ using the strong form of the Prime Number Theorem, we can use the same reasoning together with the strong form of the Prime Number Theorem for
Arithmetic Progressions to deduce that
$$
\sum^{k,\ell}\quad -\quad I^{k,\ell}<<\frac{x\log\,x\log\log\,x}{exp{(\sqrt\log\,y)}}.\tag5.23
$$
So from the above estimates, we deduce that
$$
N_2(x,k,\ell)=|S^{k,\ell}_2(x)|=\frac{x}{\phi(k)}+O(\frac{x\log\,y}{\log\,x})
+O(\frac{x\log\,x\log\log\,x}{exp{(\sqrt\log\,y)}}).\tag5.24
$$
Once again,we make the choice
$$
y=exp{\{(2\log\log\,x)^2\}},\tag5.25
$$
to deduce Theorem 7 from (5.24). 

\bigskip

\centerline{\bf{\S6: Proof of the Main Result}}
\medskip
Theorem 7 paves the way to the proof of our main result (Theorem 10 below).
Enroute to Theorem 10, we establish two theorems, the first of which
relies on Theorem 7:
\medskip
{\bf{Theorem 8:}} {\it{For integers}} $\ell, k$ {\it{satisfying}}
$1\le\ell\le k$ {\it{with}} $(\ell,k)=1$, {\it{we have}}
$$
\sum_{n\le x,\, p(n)\equiv\ell (mod\,k)}\mu(n)\omega(n)
<<\frac{x(\log\log\,x)^4}{\log\,x}.
$$
\medskip

{\bf{Proof:}} Let $f$ be a function on the primes defined by
$$
f(p)=1 \quad \text{if} \quad p\equiv\ell(mod\,k), \quad f(p)=0, \quad
\text{otherwise}. \tag6.1
$$
Then by taking $k=2$ in (1.12), and by Moebius inversion,
we get
$$
\sum_{1<n\le x, \, p(n)\equiv\ell (mod\,k)}\mu(n)(\omega(n)-1)
=\sum_{1<n\le x}\mu(n)(\omega(n)-1)f(p(n))
$$
$$
=\sum_{1<n\le x}\sum_{d|n}\mu(\frac{n}{d})f(P_2(d),
$$
which we rewrite as
$$
\sum_{1<n\le x, p(n)\equiv\ell (mod\,k)}\mu(n)\omega(n)
=\sum_{1<n\le x}\sum_{d|n}\mu(\frac{n}{d})f(P_2(d)
+\sum_{1<n\le x, \, p(n)\equiv\ell(mod\,k)}\mu(n)
$$
$$
:=\Sigma_5+\Sigma_6 \quad \text{respectively}\tag6.2.
$$
It was already established in [2] that
$$
\Sigma_6<<x\,exp\{-(\log\,x)^{(1/3)}\}.\tag6.3
$$
With regard to $\Sigma_5$, we employ the hyperbola method and write it as
$$
\Sigma_5=\sum_{m\le T}\mu(m)\sum_{d\le (x/m)}f(P_2(d))
+\sum_{d\le (x/T)}f(P_2(d))\sum_{T\le m\le (x/d)}\mu(m)
$$
$$
:=\Sigma_7+\Sigma_8, \quad \text{respectively}.\tag6.4
$$
Clearly from the strong form of the Prime Number Theorem, we get
$$
\Sigma_8<<\sum_{d\le (x/T)}f(P_2(d)\frac{x}{d\,exp{\sqrt{\log(x/d)}}}
<<\frac{x\,\log\,x}{exp{\sqrt\log\,T}}.\tag6.5
$$
Regarding $\Sigma_7$, Theorem 7 gives
$$
\Sigma_7=\sum_{m\le T}\mu(m)(\frac{x}{\phi(k)m})
+O(\frac{x(\log\log\,x)^2}{m\log(x/m)})
$$
which by the strong form of the Prime Number Theorem is
$$
<<\frac{x}{\phi(k)exp{\sqrt(\log\,T)}}+\frac{x\log\,T(\log\log\,x)^2}{\log(x/T)}.\tag6.6
$$
At this point we choose
$$
T=exp\{4(\log\log\,x)^2\} \quad <=> \quad \sqrt\log\,T=2\log\log\,x.\tag6.7
$$
With this choice of $T$, we deduce from (6.2) - (6.6) that
$$
\sum_{n\le x,\,  p(n)\equiv\ell (mod\,k)}\mu(n)\omega(n)
<<\frac{x(\log\log\,x)^4}{\log\,x}
$$
which proves Theorem 8. 

\medskip

We next prove

\medskip

{\bf{Theorem 9:}} {\it{Let}} $\ell, k$ {\it{be integers satisfying}}
$1\le \ell \le k$ {\it{with}} $(\ell, k)=1$. {\it{Then}}
$$
\sum_{1<n\le x, \, p(n)\equiv \ell(mod\,k)}\mu(n)\omega(n)\{\frac{x}{n}\}
<<\frac{x(\log\log\,x)^{5/2}}{\sqrt\,log\,x},
$$
{\it{where}} $\{y\}$ {\it{denotes the fractional part of}} $y$.  
\medskip
{\bf{Proof:}} Begin with the decomposition
$$
\sum_{1<n\le x, \, p(n)\equiv \ell(mod\,k)}\mu(n)\omega(n)\{\frac{x}{n}\}
$$
$$
=\sum_{1<n\le T, \, p(n)\equiv \ell(mod\,k)}\mu(n)\omega(n)\{\frac{x}{n}\}
+\sum_{T<n\le x, \,p(n)\equiv \ell(mod\,k)}\mu(n)\omega(n)\{\frac{x}{n}\}
$$
$$
:=\Sigma_9+\Sigma_{10} \quad \text{respectively},\tag6.8
$$
where $T$ will be chosen optimally below.

Clearly
$$
\Sigma_9<<T\log\log\,T.\tag6.9
$$
To estimate $\Sigma_{10}$, put
$$
M_{\omega}(x, \ell, k)=\sum_{1<n\le x, p(n)\equiv\ell (mod\,k)}\mu(n)\omega(n).
\tag6.10
$$
So we have
$$
\Sigma_{10}
=\sum_{T\le n\le x}(M_{\omega}(n, \ell, k)-M_{\omega}(n-1, \ell, k))\{\frac{x}{n}\}
$$
$$
=\sum_{T\le n\le x}M_{\omega}(n, \ell,k)(\{\frac{x}{n}\}-\{\frac{x}{n+1}\})
$$
$$
<<\sum_{T\le n\le x}|M_{\omega}(n,\ell,k)||\{\frac{x}{n}\}-\{\frac{x}{n+1}\}|
<<\frac{x(\log\log\,x)^4}{\log\,x}V_{\{\}}[1,\frac{x}{T}]
$$
$$
<<\frac{x(\log\log\,x)^4}{\log\,x}\frac{x}{T},\tag6.11
$$
using Theorem 8, where in (6.11), as in the proof of Theorem 2, $V_{\{\}}[a,b]$ denotes the total variation of $\{y\}$ in the interval $[a,b]$.

On comparing the bounds in (6.9) and (6.11), we see that the optimal choice of
$T$ is given by setting
$$
T\log\log\,T \sim \frac{x(\log\log\,x)^4}{\log\,x}\frac{x}{T}
$$
and this leads to the choice
$$
T=\frac{x(\log\log\,x)^{3/2}}{\sqrt\log\,x}.\tag6.12
$$
With this value of $T$, we get the upper bound in Theorem 9. 

We are now in a position to prove our main result:
\medskip
{\bf{Theorem 10:}} {\it{For integers}} $\ell, k$ {\it{satisfying}}
$1\le\ell\le k$ {\it{and}} $(\ell, k)=1$, {\it{we have}}
$$
m_{\omega}(x;\ell, k):=
\sum_{n\le x, \, p(n)\equiv\ell(mod\,k)}\frac{\mu(n)\omega(n)}{n}
<<\frac{(\log\log\,x)^{5/2}}{\sqrt\log\,x}.
$$
{\it{Letting}} $x\to\infty$, {\it{we get}}
$$
\sum_{n\ge 1, \, p(n)\equiv\ell(mod\,k)}\frac{\mu(n)\omega(n)}{n}
=\sum_{n\ge 2,\, p(n)\equiv\ell(mod\,k)}\frac{\mu(n)\omega(n)}{n}=0
$$
\medskip
{\bf{Proof:}} With $f(p)$ defined on the primes as above, note that
$$
\sum_{1<d\le x}\mu(d)(\omega(d)-1)f(p(d))[\frac{x}{d}]
=\sum_{n\le x}\sum_{1<d|n}\mu(d)(\omega(d)-1)f(p(d))
$$
$$
=\sum_{n\le x}f(P_2(n))=\frac{x}{\phi(k)}+O(\frac{x(\log\log\,x)^2}{\log\,x})
\tag6.13
$$
by (1.12) and Theorem 7.

It was already established in [2] that
$$
\sum_{1<n\le x}\mu(d)f(p(d))[\frac{x}{d}]=\sum_{n\le x}\sum_{1<d|n}\mu(d)f(p(d))
$$
$$
=-\sum_{n\le x}f(P(n))=\frac{-x}{\phi(k)}+O(\frac{x}{exp\{(\log\,x)^{1/3}\}}).
\tag6.14
$$
On comparing (6.13) and (6.14), we see that the main term $x/\phi(k)$ cancels,
and this leads to
$$
\sum_{1<d\le x}\mu(d)\omega(d)f(p(d))[\frac{x}{d}]
=O(\frac{x(\log\log\,x)^2}{\log\,x}).\tag6.15
$$
But we know by Theorem 9 that
$$
\sum_{1<n\le x}\mu(d)\omega(d)f(p(d))\{\frac{x}{d}\}
=O(\frac{x(\log\log\,x)^{5/2}}{\sqrt\log\,x}).\tag6.16
$$
Finally by adding the expressions in (6.15) and (6.16), we get
$$
x\sum_{1<n\le x}\frac{\mu(d)\omega(d)f(p(d))}{d}
=O(\frac{x(\log\log\,x)^{5/2}}{\sqrt\log\,x}).
\tag6.17
$$
On dividing both sides of (6.17) by $x$, we get Theorem 10. 

\medskip

{\bf{Remarks: The estimates when $k=1$}}
\medskip

It is to be noted that when $k=1$, the method of this section would yield quantitative estimates for $M_{\omega}(x)$, and $m_{\omega}(x)$, but these would be weaker than what we got in Theorems 1 and 4 in Section 2. First observe that when $k=1$,
the uniform distribution of $P_2(n) (mod\,1)$ is trivial, but then the error term in Theorem 7 for $k=1$ would be sharper. More precisely,
$$
N_2(x,1,1)=x+O(\frac{x}{\log\,x}).\tag6.18
$$
If we work through the proof of Theorem 8 for $k=1$ with $f$ being the characteristic function of the primes, and use (6.18), the following estimates would
hold:

The bound in (6.3) would be
$$
\Sigma_6<<x\,exp\{-c\sqrt{\log\,x}\}.\tag6.19
$$
The bound for $\Sigma_8$ in (6.5) would not change, but the bound for
$\Sigma_7$ in (6.6) would be
$$
\Sigma_7<<\frac{x}{exp{\sqrt{\log\,T}}}+\frac{x\,\log\,T}{\log(x/T)}.\tag6.20
$$
With $T$ chosen as in (6.7), the final estimate that we would get is
$$
M_{\omega}(x)<<\frac{x(\log\log\,x)^2}{\log\,x},\tag6.21
$$
which is sharper than Theorem 8 when $k=1$, but weaker that Theorem 1 which was
proved by a different method.

If we now follow the proof of Theorem 9, and use (6.21) instead of Theorem 8
for $k=1$, the optimal choice of $T$ would be
$$
T=\frac{x\,\sqrt{\log\log\,x}}{\sqrt{\log\,x}}\tag6.22
$$
in place of $T$ in (6.12). This would then yield
$$
\sum_{2\le n\le x}\mu(n)\omega(n)\{\frac{x}{n}\}<<
\frac{x(\log\log\,x)^{3/2}}{\sqrt{\log\,x}},\tag6.23
$$
which is stronger than Theorem 9 for $k=1$, but weaker than Theorem 2.  

Finally, if we use (6.23) in the proof of Theorem 10 for $k=1$, we would get
$$
m_{\omega}(x)<<\frac{(\log\log\,x)^{3/2}}{\sqrt{\log\,x}},\tag6.24
$$
which is sharper than Theorem 10 for $k=1$, but weaker than Theorem 4.

\bigskip

\centerline{\bf{\S7: Sums involving the exceptional primes}}
\medskip
In \S2 we proved (Theorem 4) that 
$$
\sum^{\infty}_{n=2}\frac{\mu(n)\omega(n)}{n}=0,\tag7.1
$$
by establishing some preliminary results. Then using the methods of \S2,
and by establishing several results, we proved in \S6 that if $k\ge 2$ is an
arbitrary modulus, then for every $\ell$ that satisfies $(\ell,k)=1$
$$
\sum^{\infty}_{n=2, \, p(n)\equiv\ell(mod\,k)}\frac{\mu(n)\omega(n)}{n}=0.\tag7.2
$$
When we sum the expression on the left in (7.2) over all $1\le\ell<k$ with
$(\ell,k)=1$, we do not get the full sum in (7.1) because the primes
$$
p\equiv\ell(mod\,k) \quad \text{with} \quad (\ell,k)>1,\tag7.3
$$
have not been accounted for. But there will be primes satisfying the conditions in (7.3),
which we call {\it{exceptional primes}}, if and only if $\ell$ is a prime divisor of $k$, and in this case there is just a single prime $p$ in the residue class
$\ell(mod\,k)$, namely $p=\ell$. It turns out that the sum in (7.2) is 0 when
taken over $n$ satisfying $p(n)=p$ for any fixed prime regardless of whether
$p$ divides $k$ or not. That is we have
\medskip
{\bf{Theorem 11:}} {\it{Let}} $p$ {\it{be an arbitrary but fixed prime. Then}}
$$
\sum^{\infty}_{n=1, \, p(n)=p}\frac{\mu(n)\omega(n)}{n}=
\sum^{\infty}_{n=2, \, p(n)=p}\frac{\mu(n)\omega(n)}{n}=0.
$$
\medskip
{\bf{Proof:}} The square-free integers $n$ with $p(n)=p$ are those of the form
$$
n=mp, \quad \text{with} \quad (m,N_p)=1, \quad  \text{where} \quad
N_p=\prod_{q\le p, \, q=\, prime}q.\tag7.4
$$
Thus using $\omega(mp)=\omega(m)+1$, we get
$$
\sum^{\infty}_{n=2, \, p(n)=p}\frac{\mu(n)\omega(n)}{n}
=\frac{-1}{p}\sum_{(m,N_p)=1}\frac{\mu(m)\omega(mp)}{m}
$$
$$
=-\frac{1}{p}\sum_{(m,N_p)=1}\frac{\mu(m)}{m}
-\frac{1}{p}\sum_{(m,N_p)=1}\frac{\mu(m)\omega(m)}{m}
$$
$$
\Sigma_{11}+\Sigma_{12}.\tag7.5
$$
It is a classical result the $\Sigma_{11}=0$. The methods of $\S2$ can be used to show that $\Sigma_{12}=0$. Thus Theorem 11 follows from (7.5).

\newpage

{\bf{Remarks:}}

(i) Since the exceptional primes, namely those that divide the modulus $k$, are finite in number, the sum of the expression in Theorem 11 taken over all exceptional primes is 0 since
it is a sum of a finite number of zeros. Thus by Theorem 11, the exceptional primes are accounted for in the full sum in (7.1). When
$(\ell, k)=1$, there are infinitely many primes $p\equiv\ell(mod\,k)$, and for each $p$ in the residue class $\ell(mod\,k)$, the sum as in Theorem 11 is 0. What makes Theorem 10 interesting is that we are summing ``infinitely many zeros'', yet the sum is 0.

(ii) As was the case with our earlier theorems, a quantitative version of
Theorem 11 can be established.

\bigskip

\centerline{\bf{\S8: The case of general $f$}}

\medskip

In the penultimate section of [2], it was shown that if $f$ is ANY bounded
function on the primes, then
$$
M_f(x):=\sum_{2\le n\le x}\mu(n)f(p(n))=o(x).\tag8.1
$$
From (8.1), it follows by Axer's theorem that
$$
\sum_{2\le n\le x}\mu(n)f(p(n))\{\frac{x}{n}\}=o(x),\tag8.2
$$
where $\{t\}$ denotes the fractional part of $t$. Next, by the Duality identity
(1.3) we have
$$
\sum_{2\le n\le x}\mu(n)f(p(n))[\frac{x}{n}]=-\sum_{2\le n\le x}f(P(n)).\tag8.3
$$
Hence by adding the expressions in (8.2) and (8,3), we get
$$
x\sum_{2\le n\le x}\frac{\mu(n)f(p(n)}{n}=-\sum_{2\le n\le x}f(P(n)) + o(x).\tag8.4
$$
From (8.4), the equivalence of (1.5) and (1.6) follows, and this was how this
equivalence was proved in [2].

In [2], the following simple bound
$$
M_f(x)<<\frac{x}{\log\log\log\,x}\tag8.5
$$
was established, but subsequently in [3] it was refined to
$$
\max_{|f|\le 1|}{|M_f(x)|}\sim\frac{2x}{\log\,x}.\tag8.6
$$
Of course, for specific functions $f$, such as $f$ being the characteristic
function of primes in an arithmetic progression $\ell(mod\,k)$, where
$(\ell, k)=1$, the bound for $M_f(x)$ is vastly superior (see [2]).

Similar in spirit to (8.1), it can be shown that

\medskip

{\bf{Theorem 12:}} {\it{If}} $f$ {\it{is any bounded function on the primes, then}}
$$
M_{f,\omega}(x):=\sum_{n\le x}\mu(n)\omega(n)f(p(n))=o(x).
$$

\medskip

A proof of a quantitative version of Theorem 12 will be given my paper with
Alamoudi [1]. From Theorem 12, by Axer's Theorem arguments, it will follow that
$$
\sum_{n\le }\mu(n)\omega(n)f(p(n))\{\frac{x}{n}\}=o(x).\tag8.7
$$
While all this seems to be similar to (8.1) and (8,2), an important difference
occurs here. In order to apply the Duality identity (1.12) when $k=2$, we have
to consider the sum
$$
\sum_{2\le n\le x}\mu(n)(\omega(n)-1)f(p(n))[\frac{x}{n}]=
\sum_{2\le n\le x} \sum_{1<d|n}\mu(d)(\omega(d)-1)f(p)d)
=\sum_{2\le n\le x}f(P_2(n)).\tag8.8
$$
We note that the first sum on the left hand side of (8.8) is
$$
\sum_{2\le n\le x}\mu(n)\omega(n)f(p(n))[\frac{x}{n}]
-\sum_{2\le n\le x}\mu(n)f(p(n))[\frac{x}{n}]
$$
$$
=\sum_{2\le n\le x}\mu(n)\omega(n)f(p(n))[\frac{x}{n}]+\sum_{2\le n\;e x}f(P(n))
\tag8.9
$$
in view of (8.3). So from (8.7), (8.8), and (8.9), we get
$$
x\sum_{2\le n\le x}\frac{\mu(n)\omega(n)f(p(n))}{n}=
\sum_{2\le n\le x}f(P_2(n)) -\sum_{2\le n\le x}f(P(n)).\tag8.10
$$
Now (8.10) yields the following result:
\medskip
{\bf{Theorem 13:}} {\it{If}} $f$ {\it{is a bounded function on the primes such that}}
$$
\sum_{2\le n\le x}f(P(n))\sim \kappa x\tag8.11
$$
{\it{and}}
$$
\sum_{2\le n\le x}f(P_2(n))\sim \kappa x,\tag8.12
$$
{\it{for some constant}} $\kappa$, {\it{then}}
$$
\sum^{\infty}_{n=2}\frac{\mu(n)\omega(n)f(p(n))}{n}=0.\tag8.13
$$
\medskip
{\bf{Remarks:}} When $f(p)$ defined on primes $p$ is the characteristic function of primes in the residue class $\ell(mod\,k)$, were $(\ell, k)=1$, then
$\kappa =1/{\phi(k)}$ in (8.11) and (8.12), in which case (8.13) is Theorem 10.  If instead of the same contant $\kappa$ in (8.12) and (8.13), we had
two different constants, $\kappa_1$ in (8.11) and $\kappa_2$ in (8.12), then
the sum in (8.13) will converge to $\kappa_2 - \kappa_1$. But we wish to stress
that we know of no natural example of a bounded function on the primes for which the constants $\kappa_1$ and $\kappa_2$ have different values. Thus we pose 
\medskip
{\underbar{PROBLEM:}} {\it{Does there exist a bounded function f on the primes such that}}
$$
\sum_{2\le n\le x}f(P(n)) \sim \kappa_1x, \quad \text{and} \quad 
\sum_{2\le n\le x}f(P_2(n)) \sim \kappa_2x, \quad \text{with} \quad
\kappa_1\neq \kappa_2.
$$

Consider now the following situation: Given $x$
arbitrarily large, define a function $f$ on the primes as follows:
$$
f(p)=1 \,\,\text{if} \,\, \sqrt x<p\le x, \quad f(p)=0 \,\, \text{if} \,\,
p\le \sqrt x.\tag8.14
$$
With $f$ as in (8.14), we have
$$
\sum_{2\le n\le x}f(P(n))=\sum_{\sqrt x<p\le x}\sum_{n\le x, P(n)=p}1
$$
$$
\sum_{\sqrt x<p\le x}[\frac{x}{p}]=xlog\,2+O(\frac{x}{log\,x}).
\tag8.15
$$
On the other hand, since $P_2(n)\le \sqrt x$ if $n\le x$, we clearly have
$$
\sum_{2\le n\le x}f(P_2(n))=0.
$$
So in this example, $\kappa_1=log\,2$, and $\kappa_2=0$. But note that the
definition of $f$ in (8.14) depends on $x$, whereas in Problem 1 we ask for a
function $f$ just defined on the primes (without dependency on $x$).

The importance of the consideration of general functions $f$ in this section will
be clear in the next section when we will discuss algebraic analogues to
the results of Alladi [2] by various authors, and algebraic analogues of the results in this paper by Sengupta [15]. 

\newpage

\centerline{\bf{\S9: Algebraic and $q$-analogues, and higher order duality}}

\medskip

The Duality identity (1.3), and the result (1.8) established in quantitative
form in Alladi [2], have attracted a lot of attention in the last decade. It all started with the paper [5] of Dawsey who obtained the following algebraic
analogue and extension of (1.8) to Galois extensions of the field of rationals
$\Cal Q$: 

Let $K$, be a Galois extension of $\Cal Q$, and ${\Cal O}_K$ the ring of integers in $K$. If $p$ is a prime in the integers, then let $P$ denote the prime
ideal that is contained in ${\Cal O}_K$ which lies above $p$. If $p$ is
unramified, let $[\frac{K/{\Cal Q}}{P}]$ denote the Artin symbol. For
simplicity, let
$$
[\frac{K/{\Cal Q}}{p}]:=[\frac{K/{\Cal Q}}{P}].
$$
Then

\medskip

{\bf{Theorem D}} (Dawsey): {\it{Let}} $K$ {\it{be a finite Galois extension
of}} $\Cal Q$ {\it{with Galois group}} $G=Gal(K/{\Cal Q})$. {\it{Let}} $C$
{\it{be a conjugacy class in}} $G$. {\it{Then}}
$$
-\sum_{n\ge 2, [\frac{K/{\Cal Q}}{p(n)}]=C}\frac{\mu(n)}{n}=\frac{|C|}{|G|}.
$$

Dawsey notes that Theorem D is a generalization of (1.8), because in the special case when $K$ is a cyclotomic extension of $\Cal Q$, the group $Gal(K/{\Cal Q})$
can be identified with $\Cal Z^*_k$, the set of reduced residues modulo $k$ for some positive integer $k$; since $\Cal Z^*_k$ is Abelian, each conjugacy class has
just one element and so
$$
\frac{|C|}{|G|}=\frac{1}{\phi(k)}.
$$

The way Dawsey proves Theorem D is to use the Chebotarev Density Theorem to show that

$$
\sum_{2\le n\le x, [\frac{K/{\Cal Q}}{P(n)}]=C}1\sim \frac{|C|}{|G|}.\tag9.1
$$
With (9.1) established, then by the methods in [2] that involve Duality, Dawsey
is able to get Theorem D.

Motivated by Dawsey's work, Sweeting and Woo [16] obtained a generalization of Theorem D in which the finite extensions $K$ of $\Cal Q$ are replaced by finite extensions $L$ of an arbitrary albegraic number field $K$. In discussing this more general situation, Sweeting and Woo consider a generalization of the Moebius function defined in terms of products of prime ideals instead of product of primes,
and establish a duality identity that generalizes (1.3) appropriately. In this more general situation, the Chebotarev Density Theorem applies, and so an analogue of Theorem D is established in [5]. 

While it is true that Theorem D generalizes (1.8), it is to be noted that the
more general equivalence of (1.5) and (1.6) is established as Theorem 6 in [2].
So what Dawsey confirmed is that if $f$ is chosen to be the characteristic
function of primes $p$ for which the Artin symbol $[\frac{K/{\Cal Q}}{P}]=C$,
then the average of $f(P(n))$ exists. That is, in this case $c$ in (1.5) is
$|C|/|G|$. So the deduction of Theorem D from (9.1) is a special case of the
equivalence of (1.5) and (1.6). Since, the equivalence of (1.5) and (1.6) is extablished in [2] for arbitrary bounded functions $f$, the bounds for the quantitative version of (1.6) is weak. For the Chebotarev Density Theorem, Lagarias
and Odlyzko [11] have established a strong form, with the error term comparable
to the error term in the strong form of the Prime Number Theorem. Thus utilizing the Lagarias-Odlyzko theorem, Dawsey is able to get a superior quantitative
version of Theorem D where the bound is just as sharp as the quantitative
version of (1.8) that is proved in [2] using the strong form of the Prime
Number Theorem.

The results of Sweeting and Woo have been extended by Kural, McDonald and Sah
[10]. A generalization in a different direction, namely replacing the Moebius
function by the more general Ramanujan sum
$$
c_m(n)=\sum^n_{k=1, (k,n)=1}e^{2imk{\pi}/n},\tag9.2
$$
is considered by Wang [21] ($\mu(n)=c_1(n)$). Also, Wang in collaboration with Duan and Yi [22] has
discussed analogues of Alladi's duality in global function fields.

A fruitful way to generalize arithmetic results is to obtain suitable $q$-analogues. In two papers [12] and [13], Ono-Schneider-Wagner have discussed a variety of
$q$-analogues of arithmetic density results and their partition
implications. 

With regard to the arithmetic consequences of the second order duality (namely
consequences of (1.12) in the case $k=2$), recently Sengupta [15], motivated by
the work of Dawsey, has obtained the extension of Theorem 10 to the situtation
when $K$ is a finite Galois extension of $\Cal Q$. Like Dawsey, Sengupta
uses the strong form of the Chebotarev Density Theorem due to Lagarias and
Odlyzko [11].

We mention that Alladi and Sengupta [4] have very recently considered arithmetic
consequences of higher order dualities, namely (1.12) for $k\ge 3$, and established analogues of all the results in this paper for $k\ge 3$. In this discussion
of higher order dualities, it turns out that when $k\ge 3$, the bounds for certain terms have extra factors which are powers of $\log\alpha$, where $\alpha=\log\,x/\log\,y$;
these factors are not present in the case $k=2$ treated here. 

Finally, we point out that all the quantitative results in [2] were established
with uniformity for the moduli $k$ of arithmetic progressions satisfying
$k\le \log^{\beta}x$, with implicit contants depending on $\beta$. This is because in [2], we utilized the Siegel-Walfisz theorem for primes in arithmetic
progressions. If we had used the Siegel-Walfisz theorem here, then Theorem 10
would hold with uniformity for $k\le \log^{\beta}x$.

\bigskip

{\bf{Concluding Remarks}:}
\medskip
(i) {\it{Arithmetic density versions:}}
\smallskip
The generalizations of (1.8) to algebraic number fields by various authors starting with Dawsey [] was motivated by rewriting (1.8) as

$$
-\sum_{n\ge 2, \, p(n)\equiv\ell(mod\,k)}\frac{\mu(n)}{n}=\frac{1}{\phi(k)},
\tag9.3
$$
and interpreting this as an arithmetic density result. Similarly, our Theorem 10 can be rewritten as
$$
\sum_{n\ge 2, p(n)\equiv\ell(mod\,k)}\frac{\mu(n)(\omega(n)-1)}{n}
=\frac{1}{\phi(k)},\tag9.4
$$
and interpreted as an arithmetic density result, thereby lending itself to an
arithmetic density generalization to algebraic number fields using the
Chebotarev density theorem (see Sengupta [15]). The consequence of the general
identity (1.12) for $k\ge 3$ discussed in
Alladi-Sengupta [4] also has an arithmetic density formulation, namely
$$
(-1)^k\sum_{n\ge 2, p(n)\equiv\ell(mod\,j)}\frac{\mu(n)}{n}\binom{\omega(n)-1}{k-1}
=\frac{1}{\phi(j)}.\tag9.5
$$
This can be generalized to algebraic number fields using the Chebotarev density theorem.

\medskip

{\it{Tenenbaum's generalization of Theorem 10:}}
\smallskip

Very recently, Tenenbaum [20] has generalized Theorem 10 as follows:

{\it{Theorem T:}} {\it{Let}} $\Cal P$ {\it{be a set of primes satisfying}}
$$
\varepsilon(t)=\frac{1}{t}\{\sum_{p\le t, \, p\in\Cal P}\log\,p\} \, -\kappa=o(1),
\quad \text{as} \quad t\to\infty,\tag9.6
$$
{\it{with some}} $\kappa\in [0,1]$. {\it{Then}}
$$
\sum^{\infty}_{n\ge 2, \, p(n)\in\Cal P}\frac{\mu(n)\omega(n)}{n}=0.\tag9.7
$$
Tenenbaum's proof of a quantitative form of (9.7) is analytic and quite intricate. But the main thing is that he is able to get (9.7) directly from (9.6)
without relying on estimates like (8.11) and (8.12). But then, our approach
using Duality connecting sums involving $\mu(n)\omega(n)f(p(n))$ with
$f(P_1(n))$ and $f(P_2(n))$ is of intrinsic interest, and that is the motivation of the present paper. 

\medskip

{\bf{Acknowledgements:}} KA would like to thank Gerald Tenenbaum for several
helpful suggestions and critical comments.

\newpage

\centerline{\bf{References}}
\medskip
1) Y. ALAMOUDI and K. ALLADI, ``Asymptotic estimates for sums involving the Moebius function,
and the number of prime factors, with conditions on the smallest prime factor'' 
(in preparation)
\smallskip
2) K. ALLADI, ``Duality between prime factors and an application to the Prime
Number Theorem for Arithmetic Progressions'', {\it{J. Num. Th.}}, {\bf{9}}
(1977), 436-451.
\smallskip
3) K. ALLADI, ``Asymptotic estimates of sums involving the Moebius function'',
{\it{J. Num. Th.}}, {\bf{14}} (1982), 86-98.
\smallskip
4) K. ALLADI and S. SENGUPTA, ``Higher order duality between prime factors
and primes in arithmetic progressions'', (in preparation).
\smallskip
5) M. L. DAWSEY, ``A new formula for Chebotarev densities'', {\it{Res. Num. Th.}}, {\bf{3}} (2017).

https://doi.org/10.1007/s40993-017-0093-7

\smallskip
6) N. G. DE BRUIJN, ``On the asymptotic behavior of a function occurring in
the theory of primes'', {\it{J. Indian Math. Soc.} (N.S.)}, {\bf{15}} (1951),
25-32.
\smallskip
7) N. G. DE BRUIJN, ``On the number of positive integers $\le x$ and free of
prime factors $>y$'', {\it{Indag. Math.}}, {\bf{13}} (1951), 50-60.
\smallskip
8) L. DUAN, B. WANG, and S. YI, ``Analogues of Alladi's formula over global
function fields'', {\it{Finite Fields and their Applications}}, {\bf{74}}
(2021),

https://doi.org/10.1016/j.ffa.2021.101874
\smallskip
9) A. J. HILDEBRAND and G. TENENBAUM, ``On integers free of large prime factors'', {\it{Trans. Amer. Math. Soc.}}, {\bf{296}} (1986), 265-290.
\smallskip
10) M. KURAL, V. MCDONALD, and A. SAH, ``Moebius formulas for densities of
sets of prime ideals'', {\it{Arch. Math.}} {\bf{135}} (2020), 53-60.
\smallskip
11) J. C. LAGARIAS and A. M. ODLYZKO, ``Effective versions of the Chebotarev
density theorem'', in {\it{Algebraic Number Fields, $L$-functions, and Galois
properties}}, A. Fr\"ohlich (Ed.), Acad. Press. London (1977), 409-464.
\smallskip
12) K. ONO, R. SCHNEIDER, and I. WAGNER, ``Partition theoretic formulas for
arithmetic densities'', in {\it{Analytic Number Theory, Modular Forms, and $q$-Hypergeometric Series}} (G. E. Andrews and F. Garvan, Eds.) - Conf. in honor of
Krishna Alladi's 60th birthday, Springer Proceedings in Math. and Stat., {\bf{221}} (2017), 611-624.
\smallskip
13) K. ONO, R. SCHNEIDER, and I. WAGNER, ``Partition-theoretic formulas for
arithmetic densities - II'', {\it{Hardy-Ramanujan J.}}, {\bf{43}} (2020), 1-16.
\smallskip
14) A. SELBERG, ``Note on a paper by L. G. Sathe'', {\it{J. Indian Math. Soc.}},
{\bf{18}} (1954), 83-87.
\smallskip
15) S. SENGUPTA, ``Algebraic analogues of theorems of Alladi-Johnson relating
to second order duality among prime factors'' (in preparation)
\smallskip
16) N. SWEETING and K. WOO, ``Formulas for Chebotarev densities of Galois
extensions of number fields'', {\it{Research Num. Th.}} {\bf{5}} (2019) -

https://doi.org/10.1007/s40993-018-0142-x
\smallskip
17) G. TENENBAUM, ``A rate estimate in Billingsley's theorem for the size
distribution of the number of prime factors'', {\it{Quart. J. Math.}},
{\bf{51}}
(2000), 385-403.
\smallskip
18)  G. TENENBAUM, {\it{Probabilistic and Analytic Number Theory}}, Grad. Studies in Math., {\bf{163}} (2015), Amer. Math. Soc, Providence RI, 630 pp. 
\smallskip
19) G. TENENBAUM, ``Private Communication to Alladi'' (2019)
\smallskip
20) G. TENENBAUM, ``On a family of arithmetic series related to the Moebius function'' (2024 preprint).
\smallskip
21) B. WANG, ``The Ramanujan sum and Chebotarev densities'', {\it{Ramanujan J.}}, {\bf{55}} (2021), 1105-1111.
\smallskip
22) B. WANG, ``Analogues of Alladi's formula'', {\it{J. Num. Th.}} {\bf{221}}
(2021), 232-246.

\bigskip

Department of Mathematics

University of Florida

Gainesville, FL 32611

USA

\medskip

email:

KRISHNASWAMI ALLADI - alladik(at)ufl.edu

JASON JOHNSON - iridiumalchemist(at)gmail.com

\end